\documentclass[draft,12pt, openany, a4paper] {article}
\usepackage{amsmath}
\usepackage{amssymb}

\begin{document}

\title{ \Large {The Complex Structures on  $S^{2n}$}
 \footnotetext{
\mbox{}\quad  \ \ {\it Subject classification}. \  53B21, 53B35,   53C28,57R20.\\
\mbox{}\quad  \ \ E-mail: \ jwzhou@suda.edu.cn }}
\author{Jianwei Zhou}
\date{\small Department of Mathematics, Suzhou University, Suzhou
215006, P. R. China }

\maketitle
\begin{abstract} Let $\widetilde{\cal
J}(S^{2n})$ be the set of orthogonal  complex structures on
 $TS^{2n}$. We show that the twistor space $\widetilde{\cal
J}(S^{2n})$  is a Kaehler manifold. Then we show that an orthogonal
almost complex structure $J_f$ on $S^{2n}$ is integrable if and only
if the corresponding section $f\colon\; S^{2n}\to \widetilde{\cal
J}(S^{2n}) $ is holomorphic. These shows there is no integrable
orthogonal complex structure on the sphere $S^{2n}$ for $n>1$. We
also show that there is no  complex structure in a neighborhood of
the space $\widetilde{\cal J}(S^{2n})$.  The method is to study the
first Chern class of $T^{(1,0)}S^{2n}$.
\end{abstract}

\baselineskip 15pt
\parskip 3pt

\noindent{\bf  1. Introduction } \vskip 0.2cm

An almost complex structure $J$ on a differentiable manifold is an
endomorphism  of the tangent bundle $TM$ such that $J^2=-1$. It is
known [2] that the only sphere admitting such structures are $S^2$
and $S^6$, see also Theorem 10.4.7 of [11] where  the index theorem
is used.

The operator $J\colon\; TM\to TM$ is said to be integrable if it
comes from an honest complex structure on $M$, then $M$ is a complex
manifold. It is a long-standing problem that whether there is an
integrable complex structure on the sphere $S^{6}$.  In 1987 Lebrun
[4]  proved $\nabla_{X_\alpha}X_\beta = - \nabla_{X_\beta}X_\alpha$
for $(0,1)$ vector fields  if the orthogonal complex structure $J$
is integrable. Then he claimed that this can be used to prove that
the map $\tau\colon\; U \subset S^6\to G_C(3,7)$ is a holomorphic
map. But in \S 4 we show that $\nabla_{X_\alpha}X_\beta =-
\nabla_{X_\beta}X_\alpha$ can not hold.

In [7], we use  Clifford algebra and the spinor calculus to study
the orthogonal complex structures on  the sphere $S^6$. By the spin
representation  we show that the Grassmann manifold $G(2,8)$ can be
viewed as the set of orthogonal complex structures on $R^8$ and
there is a fibre bundle $\tau \colon\; G(2,8)\to S^{6}$.  Then we
show that there is  no orthogonal complex structure on the sphere
$S^6$.

In this paper we study the general problem. Let ${\cal J}(R^{2n})$
be the set of all complex structures on Euclidean space $R^{2n}$. In
\S 2, we show that with naturally defined metric and complex
structure, the twistor space ${\cal J}(R^{2n})$ is a Kaehler
manifold and the orthogonal twistor space $\widetilde{\cal
J}(S^{2n})=\widetilde{\cal J}(R^{2n+2})$ is a Kaehler submanifold of
${\cal J}(R^{2n+2})$.

There are two connected components on ${\cal J}(R^{2n})$. Denote
${\cal J}^+(R^{2n})$ the space of all oriented complex structures
and we have $\widetilde{\cal J}^+(S^{2n})=\widetilde{\cal
J}^+(R^{2n+2})$. Theorem 2.6 shows that the Poincar\'{e} polynomial
of $\widetilde{\cal J}^+(S^{2n})$ is
$$P_t(\widetilde{\cal J}^+(S^{2n}))
=(1+t^2)(1+t^4)\cdots (1+t^{2n}).$$

In \S 3, we show that an orthogonal almost complex structure $J$ on
$S^{2n}$ is integrable if and only if the corresponding section
$f\colon\; S^{2n}\to \widetilde{\cal J}(S^{2n})$ is holomorphic.
These shows there is no orthogonal complex structure on the sphere
$S^{2n}$ for $n>1$.

In \S 4, we use Riemannian connection on $S^{2n}$ to study the
problem, show that a local section $f\colon\; S^{2n}\to {\cal
J}(S^{2n})\subset {\cal J}(R^{2n+2})$ is holomorphic if and only if
the complex structure defined by $f$ is orthogonal and integrable.

In \S 5, we study the map $f: S^{2n}\to G_{\bf C}(n, 2n+1)$ defined
by Lebrun. In \S 6, we show that there is no  complex structure in a
neighborhood of the space $\widetilde{\cal J}(S^{2n}) (n>1)$.  The
method is to study the first Chern class of vetcor bundle
$T^{(1,0)}S^{2n}$.

\vskip 0.7cm \noindent{\bf  2. The twistor space ${\cal J}(R^{2n})$}
\vskip 0.2cm

Let $GL(2n,R)$ be the general linear group and $J_0=\left(
\begin{array}{cccccccc}
0  & -1  & &  &  &   \\
1  &0  &  &  &  &    \\
    &   &    \ddots  &   &   \\
   &  &     &   &0  & -1  \\
  &  &     &   & 1  &0  \\
\end{array}\right)$
be a complex structure on  Euclidean space $R^{2n}$. Let $GL(n,C)$
be the stability subgroup of $J_0$:
$$GL(n,C)=\{ G\in GL(2n,R) \ | \ GJ_0=J_0G\}.$$
Then
$${\cal J}(R^{2n})=\frac {GL(2n,R)}{GL(n,C)}=\{GJ_0G^{-1} \ | \ G\in GL(2n,R)\}$$
is the set of all complex structures on $R^{2n}$ called a twistor
space on $R^{2n}$. It is easy to see
$${\cal J}(R^{2n})=\{A\in GL(2n,R) \ | \ A^2=-I\}.$$

 {\bf Lemma 2.1}  \ For any $A=GJ_0G^{-1}\in {\cal J}(R^{2n})$, \
 $G=(e_{1},e_2,\cdots,e_{2n})$,  we have
 $Ae_{2i-1}=e_{2i}, \ \ Ae_{2i}=-e_{2i-1}$, where $e_l\in
 R^{2n}$ are column vectors.

{\bf Proof} \ The Lemma follows from
$$A(e_{1},e_2,\cdots,e_{2n})=(e_{1},e_2,\cdots,e_{2n})J_0=
(e_2,-e_1,\cdots,e_{2n},-e_{2n-1}).  \ \ \ \Box $$

For any $B\in GL(2n,R)$, define  inner product on
$T_BGL(2n,R)=gl(2n,R)$ by
$$\langle X,Y \rangle = \frac 12 \mbox{tr} \,  (XY^t)+\frac 12 \mbox{tr} \,
(BXB^{-1}(BYB^{-1})^t), \ \ X,Y\in
gl(2n,R).$$  Restricting this inner product  on $T{\cal
J}(R^{2n})$ makes ${\cal J}(R^{2n})$ a Riemannian manifold.

It is easy to see that
$$T_A{\cal J}(R^{2n}) = \{ X\in gl(2n,R) \ | \ X=AXA \ \mbox{or} \ AX=-XA\}.$$ The normal
space of $T_A{\cal J}(R^{2n})$ in $T_AGL(2n,R)$ with inner product
$\langle \ ,\rangle $ is
$$T^\perp _A{\cal J}(R^{2n})  =  \{ Y\in gl(2n,R) \ | \ Y=-AYA \ \mbox{or} \ AY=YA\}.$$ The spaces $T_A{\cal
J}(R^{2n})$  and $T^\perp _A{\cal J}(R^{2n})$ are all invariant by
the action $Ad(A)$. Then for any $X_1,X_2\in T_A{\cal J}(R^{2n})$,
we have
$$\langle X_1,X_2\rangle =  \mbox{tr} \,  (X_1X_2^t).$$
Similarly, for any $Y_1,Y_2\in T_A^\perp{\cal J}(R^{2n})$, we have
$\langle Y_1,Y_2\rangle =  \mbox{tr} \, (Y_1Y_2^t)$.

For any $A\in {\cal J}(R^{2n})$, the map $X\in gl(2n,R)\to  AX$
defines a complex structure on $gl(2n,R)$. It is easy to see that
for any $X\in T_A{\cal J}(R^{2n}), \ AX\in T_A{\cal J}(R^{2n})$.

 {\bf Lemma 2.2}  \ The maps $\tilde JX=AX (X\in T_A{\cal J}(R^{2n}))$
 define an almost complex structure $\tilde J$ on ${\cal
J}(R^{2n})$.

Let $\widetilde D$ be the  Riemannian connection on $ {\cal
J}(R^{2n})$ determined by  the   metric $\langle \ , \rangle.$ In
general, $\langle X_1,X_2\rangle \neq \langle \tilde JX_1,\tilde
JX_2\rangle, \ \ X_1,X_2\in T_A{\cal J}(R^{2n}).$ For any
$X_1,X_2\in T{\cal J}(R^{2n})$ define
$$ds^2(X_1,X_2)= \frac 12 \langle X_1,X_2\rangle  + \frac 12 \langle \tilde JX_1,\tilde JX_2\rangle .$$
The almost complex structure $\tilde J$ preserves  the metric
$ds^2$.

 {\bf Lemma 2.3}  \ $\widetilde D\tilde J=0$.

{\bf Proof}  \ Let $X_1,\cdots,X_{2n^2}$ be tangent frame fields on
${\cal J}(R^{2n})$ such that $\tilde JX_{2i-1}=X_{2i}, \ \tilde
JX_{2i}=-X_{2i-1}.$ Let $Y_\alpha$ be frame fields  on $T^\perp
{\cal J}(R^{2n})$. Then we can write
$$dX_l=\sum \ \omega_l^kX_k +\sum \ \omega_l^\alpha Y_\alpha.$$
By Gauss formula, the connection $\widetilde D$ is defined by
$$\widetilde DX_l=\sum \ \omega_l^kX_k.$$

Let $dA=\sum \ \omega^lX_l$. From $\tilde
JX_{2i-1}=AX_{2i-1}=X_{2i}$, we have
\begin{eqnarray*}dX_{2i} &=& \sum \ \omega_{2i}^kX_k +\sum \ \omega_{2i}^\alpha Y_\alpha \\
&=& \sum \ \omega^lX_lX_{2i-1}+\sum \ \omega_{2i-1}^kAX_k +\sum \
\omega_{2i-1}^\alpha AY_\alpha.
\end{eqnarray*}
By $AX_lX_{2i-1}=X_lX_{2i-1}A$, we know that $X_lX_{2i-1}\in
T^\perp {\cal J}(R^{2n})$, hence $\sum \ \omega_{2i}^kX_k=\sum \
\omega_{2i-1}^kAX_k$ and we have
$$\omega_{2i-1}^{2j-1}=\omega_{2i}^{2j}, \ \ \omega_{2i}^{2j-1}=-\omega_{2i-1}^{2j}.$$
Then $\widetilde D(\tilde JX_{2i-1})=\tilde J\widetilde DX_{2i-1}$
and
$$(\widetilde D\tilde J)X_{2i-1}= \widetilde D(\tilde JX_{2i-1})-\tilde J\widetilde DX_{2i-1}=0.$$
Similarly, $(\widetilde D\tilde J)X_{2i}= 0.$  \ \ \ $\Box$

{\bf Theorem  2.4}  \ $\tilde J$ is an orthogonal  complex structure
on ${\cal J}(R^{2n})$ with metric $ds^2$. Then  the twistor space
$({\cal J}(R^{2n}),ds^2,\tilde J)$ is a Kaehler manifold.

 {\bf Proof}  \ With the almost complex structure $\tilde J$ and the metric $ds^2$,
${\cal J}(R^{2n})$ is an almost Hermitian manifold.
 $\widetilde D$ is Riemannian connection with metric
 $\langle \ , \rangle$.
 For any $X_1,X_2\in \Gamma(T{\cal J}(R^{2n}))$, from $\widetilde D\tilde J=0$
we have $\widetilde D(\tilde JX_i)=\tilde J\widetilde DX_i$,
\begin{eqnarray*} && 2d[ds^2(X_1,X_2)] \\
&=& \langle \widetilde DX_1,X_2\rangle+ \langle X_1,\widetilde
DX_2\rangle+ \langle \widetilde D(\tilde JX_1), \tilde JX_2\rangle+
\langle \tilde J
X_1,\widetilde D(\tilde JX_2)\rangle \\
&=& \langle \widetilde DX_1,X_2\rangle+ \langle X_1,\widetilde
DX_2\rangle+ \langle \tilde J\widetilde DX_1,\tilde JX_2\rangle+
\langle \tilde JX_1,\tilde J\widetilde DX_2\rangle \\
& = &2ds^2(\widetilde DX_1,X_2)+ 2ds^2(X_1,\widetilde
DX_2).\end{eqnarray*} These shows $\widetilde D$ is  the Riemannian
connection on $({\cal J}(R^{2n}),ds^2)$. By Theorem 7.4.3 of [11],
from $\widetilde D\tilde J=0$ we know that $\tilde J$ is integrable
and the  fundamental 2-form  of $({\cal J}(R^{2n}),ds^2,\tilde J)$
is closed. Then
 ${\cal J}(R^{2n})$ is a Kaehler manifold. \
\ \ $\Box$

The orthogonal twistor space $$\widetilde {\cal J}(R^{2n})=\{GJ_0G^t
\ | \ G\in O(2n)\}=\{A\in {\cal J}(R^{2n}) \ | \ A^t=-A\}$$ is a
subspace of ${\cal J}(R^{2n})$. It is easy to see that
$$T_A\widetilde {\cal J}(R^{2n})=\{X\in gl(2n, R) \ | \ X=-X^t=AXA\}.$$
Note that on $\widetilde {\cal J}(R^{2n})$ we have $\langle \tilde
JX_1,\tilde JX_2\rangle= \langle X_1,X_2\rangle$. The induced metric
on $\widetilde {\cal J}(R^{2n})$ is
$$ds^2(X_1,X_2)=\langle X_1,X_2\rangle=\mbox {tr}\, (X_1X_2^t).$$

{\bf Theorem 2.5}  \ (1) $\widetilde {\cal J}(R^{2n})$ is a Kaehler
submanifold of ${\cal J}(R^{2n})$;

(2) $\widetilde {\cal J}(R^{2n})$ is a deformation retract of
${\cal J}(R^{2n})$.

{\bf Proof}  \ (1) \  For any $X\in T_A\widetilde {\cal J}(R^{2n})$,
we have $\tilde JX=AX\in T_A\widetilde {\cal J}(R^{2n})$, $\tilde J$
is a complex structure on $\widetilde {\cal J}(R^{2n})$. Then
$\widetilde {\cal J}(R^{2n})$ is a complex submanifold of ${\cal
J}(R^{2n})$, hence a Kaehler submanifold.

(2) \ Let $A\in {\cal J}(R^{2n})$, set $A=A_1+A_2$ where
$A_1=-A_1^t, \ A_2=A_2^t$.  From $A^2=A_1^2+A_2^2+A_1A_2+A_2A_1=-I$
and $(A_1A_2+A_2A_1)^t=-A_1A_2-A_2A_1$, we have
$$A_1^2+A_2^2=-I, \ \ A_1A_2+A_2A_1=0.$$
If there is a vector $x$ such that $A_1x=0$, we have
$-x=A^2x=A^2_2x=A^t_2A_2x$, then $-x^tx=x^tA^t_2A_2x$ and $x=0$.
Then $A_1$ is a non-singular matrix. Let $e_1,\cdots, e_{2n}$ be
orthonormal vectors such that $A_2 e_i=\lambda_i e_i$, then $\sum \
e_i\cdot e_i^t=I$ and  $A_2=\sum \ \lambda_i e_i\cdot
 e_i^t$.  Let $A_1= BP$ be the unique decomposition such that $B\in
\widetilde {\cal J}(R^{2n}), \ BP=PB, \ P^2=-A_1^2$ and $P$ is
positive definite, for proof see [6], p.205,213. Hence $P=\sum \
\sqrt {1+\lambda_i^2} \,  e_i\cdot  e_i^t, \ PA_2=A_2P, \
BA_2+A_2B=0$. The element $A$ is determined by $B$ and $A_2$ with
$BA_2+A_2B=0$.

As $A_2B e_i=-BA_2 e_i=-\lambda_i e_i$, then if $\lambda_i$ is a
non-zero characteristic value of $A_2$, $-\lambda_i$ is also a
characteristic value of $A_2$ with characteristic vector $B e_i, \
B(B e_i)=- e_i$. These shows that for any $A\in {\cal J}(R^{2n})$,
there are orthonormal vectors $e_i, e_{n+i}$ and numbers
$\lambda_i$ such that $A=A_1+A_2=BP+A_2$ with
$$A_2=\sum\limits_{i=1}^n \  \lambda_i(e_ie_i^t - e_{n+i}e_{n+i}^t),
\ \ P=\sum\limits_{i=1}^n \ \sqrt {1+\lambda_i^2} (e_ie_i^t +
e_{n+i}e_{n+i}^t),$$ $$B= \sum\limits_{i=1}^n \
 (e_{n+i}e_i^t - e_{i}e_{n+i}^t),$$
$$A=\sum\limits_{i=1}^n \ \sqrt{1+\lambda_i^2}
 (e_{n+i}e_i^t - e_{i}e_{n+i}^t) + \sum\limits_{i=1}^n \ \lambda_i (e_ie_i^t -
 e_{n+i}e_{n+i}^t).$$

Let $A_2(t)=tA_2, \ P(t)=\sum \ \sqrt {1+t^2\lambda_i^2}\,(e_ie_i^t
+ e_{n+i}e_{n+i}^t)$. Then $ A(t)= BP(t)+A_2(t)\in {\cal J}(R^{2n}),
\ A(0)=B, \ A(1)=A$. These shows  $\widetilde {\cal J}(R^{2n})$ is a
deformation retract of ${\cal J}(R^{2n}). $

Then we have a map ${\cal J}(R^{2n}) \to \widetilde {\cal
J}(R^{2n})$ defined by $A \to  B$.\ \ \ $\Box$

 Let
$e_{-1}=(1,0,\cdots,0)^t\in R^{2n+2}$ be a fixed vector and
$$S^{2n}=\{e_0\in R^{2n+2} \ | \ |e_0|=1, \ e_0\perp
e_{-1}\}.$$ Let space ${\cal J}(S^{2n})$  be the subspace of $A\in
{\cal J}(R^{2n+2})$ such that $e_0=Ae_{-1}\in S^{2n}$, \ $Av\perp
e_{-1},e_0$ for any $v\perp e_{-1},e_0$. ${\cal J}(S^{2n})$ is the
twistor space on $S^{2n}$,  the projection
 $\pi\colon\;
{\cal J}(S^{2n})\to S^{2n}$ is $\pi (A)=Ae_{-1}$.

Let $\widetilde{\cal J}(S^{2n})$ be a subspace of ${\cal J}(S^{2n})$
formed by all orthogonal complex structures on $S^{2n}$. For any
$A=GJ_0G^{-1}\in \widetilde {\cal J}(R^{2n+2}), \ G\in O(2n+2)$, we
can choose an $H\in U(n+1)$ such that $GH=(e_{-1},
e_0,e_1,\cdots,e_{2n})$. Then
$$A=GJ_0G^{t}=e_0e_{-1}^t - e_{-1}e_0^t +\sum\limits_{i=1}^n \
(e_{2i}e_{2i-1}^t - e_{2i-1}e_{2i}^t).$$ $\sum\limits_{i=1}^n \
(e_{2i}e_{2i-1}^t - e_{2i-1}e_{2i}^t)$ defines an orthogonal complex
structure on  $T_{e_0}(S^{2n})$. These shows $\widetilde {\cal
J}(S^{2n})=\widetilde {\cal J}(R^{2n+2})$, see also [3].

{\bf Theorem  2.6}  \  $(\widetilde{\cal J}(S^{2n}), d^2s, \tilde
J)$ is a Kaehler submanifold of ${\cal J}(R^{2n+2})$.

There are two connected components on $\widetilde{\cal J}(S^{2n})$.
Denote $$\widetilde{\cal J}^+(S^{2n})=\{GJ_0G^t \ | \ G\in
SO(2n+2)\}$$ the twistor space of all oriented almost orthogonal
complex structure on $S^{2n}$. Using degenerate Morse functions we
can show

{\bf Theorem 2.7}  \ The Poincar\'{e} polynomial of $\widetilde{\cal
J}^+(S^{2n})$ is
$$P_t(\widetilde{\cal J}^+(S^{2n}))=(1+t^2)(1+t^4)\cdots (1+t^{2n}).$$

{\bf Proof} \ As in [8],[9], let $$h(A)=\langle A,\bar e_0 e_{-1}^t-
e_{-1}\bar e_0^t\rangle$$ be a function on $\widetilde{\cal
J}^+(S^{2n})$, where $e_{-1}, \bar e_{0}=(0,1,0,\cdots,0)^t$ are two
fixed vectors. Let $A=e_0e_{-1}^t - e_{-1}e_0^t +\sum\limits_{i=1}^n
\ (e_{2i}e_{2i-1}^t - e_{2i-1}e_{2i}^t)$ and
$de_0=\sum\limits_{i=1}^{2n} \ \omega^ie_i.$ Hence
\begin{eqnarray*} dh & = &  \sum\limits_{i=1}^{n} \ \omega^{2i}
\langle e_{2i} e^t_{-1} - e_{-1}e^t_{2i},
\bar e_0 e_{-1}^t- e_{-1}\bar e_0^t\rangle \\
& & + \sum\limits_{i=1}^{n} \ \omega^{2i-1}\langle e_{2i-1} e^t_{-1}
-e_{-1}e^t_{2i-1},\bar e_0
e_{-1}^t- e_{-1}\bar e_0^t\rangle \\
&=& 2\sum\limits_{l=1}^{2n} \ \omega^{l}\langle e_l,\bar e_0\rangle,
\end{eqnarray*}
$dh=0$ if and only if $e_0=\bar e_0$ or $e_0=-\bar e_0$. Then
$\pi^{-1}(\bar e_0)$ and $\pi^{-1}(-\bar e_0)$ are two critical
submanifolds of the function $h$. It is easy to see that $$
d^2h|_{\pi^{-1}(\bar e_0)} = - 2\sum\limits_{l=1}^{2n} \
\omega^{l}\otimes \omega^{l}, \ \ \
 d^2h|_{\pi^{-1}(-\bar e_0)} =
2\sum\limits_{l=1}^{2n} \ \omega^{l}\otimes \omega^{l}.$$ These
shows that the critical submanifolds $h^{-1}(-2)=\pi^{-1}(-\bar
e_0)$ and $h^{-1}(2)=\pi^{-1}(\bar e_0)$ are non-degenerate with
indices $0$ and $2n$ respectively. Then $h$ is a Morse function and
  the Poincar\'{e} polynomial of $\widetilde{\cal J}^+(S^{2n})$ is
$$P_t(\widetilde{\cal J}^+(S^{2n}))=(1+t^{2n})P_t(\widetilde{\cal J}^+(S^{2n-2})).  \ \ \ \Box$$

\vskip 0.7cm \noindent{\bf  3. The orthogonal complex structure
 on the sphere $S^{2n}$} \vskip
0.2cm

In this section we use twistor space $\pi\colon\;  \widetilde{\cal
J}(S^{2n})\to S^{2n}$ to study the complex structure on $S^{2n}$.
The elements of $\widetilde{\cal J}(S^{2n})$ can be represented by
$A=\sum\limits_{i={0}}^n \ (e_{2i}\cdot e_{2i-1}^t- e_{2i-1}\cdot
e_{2i}^t),$ where $e_{-1}=(1,0,\cdots,0)^t, \ e_0\in S^{2n}, \
(e_{-1},e_0,e_1,\cdots,e_{2n})\in O(2n+2)$.

 By the method of
moving frame we have
$$d(e_0,e_1,\cdots,e_{2n})=(e_0,e_1,\cdots,e_{2n})
\left( \begin{array}{cccccccc}
0  & -\omega^1  & -\omega^2 & \cdots  &  -\omega^{2n}  \\
\omega^1  & 0  & -\omega_{12} & \cdots & -\omega_{1,2n} \\
\omega^2  & \omega_{12} & 0  & \cdots  & -\omega_{2,2n}  \\
\omega^3  & \omega_{13} & \omega_{23}  & \cdots  & -\omega_{3,2n}  \\
 \vdots  & \vdots  & \vdots   & \ddots  & \vdots \\
 \omega^{2n}  & \omega_{1,2n} & \omega_{2,2n}  & \cdots  & 0  \\
\end{array}\right).$$
 By  simple computation, we have
$$ dA =\sum\limits_{i<j} \ (\alpha_{ij}^*\alpha_{ij}
 + \beta_{ij}^*\beta_{ij})  + \sum\limits_{i=1}^{n} \ (\omega^{2i-1}\widetilde X_{2i-1}
 + \omega^{2i}\widetilde X_{2i}),
$$
where
$$\alpha_{ij}=e_{2j-1}e^t_{2i-1} -e_{2j}e^t_{2i}
-e_{2i-1}e^t_{2j-1} +e_{2i}e^t_{2j},$$ $$ \beta_{ij}
=e_{2j}e^t_{2i-1} +e_{2j-1}e^t_{2i} -e_{2i}e^t_{2j-1}
-e_{2i-1}e^t_{2j},$$
 $$\widetilde X_{2i-1}=e_{2i-1} e^t_{-1} -
e_{-1}e^t_{2i-1}+e_{0}e^t_{2i} -e_{2i}e^t_{0}, $$
$$ \widetilde X_{2i}=e_{2i} e^t_{-1} -
 e_{-1}e^t_{2i}-e_{0}e^t_{2i-1} +e_{2i-1}e^t_{0},$$
$$\alpha_{ij}^*=\omega_{2i,2j-1}+\omega_{2i-1,2j}, \ \ \ \beta_{ij}^*
=\omega_{2i,2j}-\omega_{2i-1,2j-1}.$$

 {\bf Lemma 3.1}  \ $\alpha_{ij}, \beta_{ij},\widetilde  X_{2i-1},\widetilde X_{2i}$
  are local basis of $T\widetilde{\cal J}(S^{2n})$,  and
$\alpha_{ij}^*, \ \beta_{ij}^*,   \omega^{2i-1}, \omega^{2i}$ are
the dual basis. Furthermore we have
 $$\tilde J\alpha_{ij}=A\alpha_{ij}= \beta_{ij}, \ \
\tilde J\widetilde X_{2i-1}=A\widetilde X_{2i-1} =-\widetilde X_{2i}
 $$
and $\tilde J\alpha_{ij}^*= -\beta_{ij}^*, \ \ \tilde
J\omega^{2i-1}=-\omega^{2i}.$

Now let $J$ be an orthogonal almost complex structure on $S^{2n}$.
If $2n\not= 2,6$, $J$ is defined on a open subset of $S^{2n}$. Let
$e_{-1},e_0, e_1,\cdots, e_{2n}$ be local orthonormal frame fields,
$e_{-1}=(1,0,\cdots,0)^t$, \ $e_0$ be the position vector and
$e_1,\cdots, e_{2n} \in T_{e_0}S^{2n}$ with $Je_{2i-1}=e_{2i},\
i=1,\cdots,n$. We have
$$de_0=\sum\limits_{i=1}^{2n} \ \omega^ie_i, \ \ de_i=\sum\limits_{j=1}^{2n} \
\omega_i^je_j-\omega^ie_0,$$ $$d\omega^i= \sum\limits_{j=1}^{2n} \
\omega^j\wedge\omega_j^i, \ \ \omega_j^i+\omega^j_i=0. $$ $\{e_i\}$
and $\{\omega^i\}$ are dual frame fields on $S^{2n}$. The operator
$J$ also acts on the cotangent vectors by
$J\omega^{2i-1}=-\omega^{2i}, \ J\omega^{2i}=\omega^{2i-1}$;
$\varphi^i=\omega^{2i-1}+\sqrt {-1}\omega^{2i}, \ \varphi^{\bar
i}=\omega^{2i-1}-\sqrt {-1}\omega^{2i}$ are $(1,0)$  and $(0,1)$
forms respectively.   Let $T^{(1,0)}S^{2n}$ be subspace of
$TS^{2n}\otimes C$ generated by $Z_{i}=e_{2i-1}-\sqrt {-1}e_{2i}, \
\ i=1,\cdots,n$; \ $T^{(0,1)}S^{2n}$ be  generated by $Z_{\bar
i}=e_{2i-1}+\sqrt {-1}e_{2i}.$

By $d\omega^i= \sum\limits_{j=1}^{2n} \ \omega^j\wedge\omega_j^i$ we
have $d\varphi^i=\sum\limits_{j=1}^n \ \varphi^j\wedge\varphi_j^i+
\sum\limits_{j=1}^n \ \varphi^{\bar j} \wedge\varphi_{\bar j}^i,$
where {\footnotesize
$$\left( \begin{array}{cc}
\varphi_{j}^i  & \varphi_{j}^{\bar i}  \\
\varphi_{\bar j}^i & \varphi_{\bar j}^{\bar i}
\end{array}\right)=\frac 12 \left( \begin{array}{cccccccc}
\omega_{2j-1}^{2i-1}+
\omega_{2j}^{2i}+\sqrt{-1}(\omega_{2j-1}^{2i}-\omega_{2j}^{2i-1} ) &
\omega_{2j-1}^{2i-1}- \omega_{2j}^{2i}-\sqrt{-1}(\omega_{2j-1}^{2i}+\omega_{2j}^{2i-1} ) \\
\omega_{2j-1}^{2i-1}-
\omega_{2j}^{2i}+\sqrt{-1}(\omega_{2j-1}^{2i}+\omega_{2j}^{2i-1} ) &
\omega_{2j-1}^{2i-1}+
\omega_{2j}^{2i}-\sqrt{-1}(\omega_{2j-1}^{2i}-\omega_{2j}^{2i-1} )
\end{array}\right).  $$ }
It is easy to see that
$$dZ_{i}=\sum\limits_{j=1}^n \ (\varphi_{i}^{j} Z_{j}+\varphi_{i}^{\bar j} Z_{\bar j})
-\varphi^{\bar i}e_0.$$ Then from Riemannian connection on $S^{2n}$,
we have  $\nabla Z_{i}=\sum\limits_{j=1}^n \ (\varphi_{i}^{j}
Z_{j}+\varphi_{i}^{\bar j} Z_{\bar j}).$ As [3], [7], we have

 {\bf Lemma 3.2}  \ $J$ is integrable if and only if every $\nabla_{Z_j}Z_i$ is $(1,0)$ vector field,
 or equivalently, every
 $\varphi_i^{\bar j}$ is $(0,1)$-form, then $J(\omega_{2j-1}^{2i}+\omega_{2j}^{2i-1})=
\omega_{2j-1}^{2i-1}- \omega_{2j}^{2i}$.

{\bf Proof} \  The Lemma can also be proved as [1]:

Assuming $J$ is integrable, for any
$X,Y,Z\in\Gamma(T^{(1,0)}S^{2n})$, we have
$[X,Y]=\nabla_XY-\nabla_YX\in \Gamma(T^{(1,0)}S^{2n})$ and $\langle
X, Y\rangle =0$ with respect to the standard metric on $S^{2n}$.
Then $0=Z\langle X, Y\rangle =\langle \nabla_ZX, Y\rangle+\langle X,
\nabla_ZY\rangle,$ hence
\begin{eqnarray*} \langle \nabla_XY, Z\rangle & = & \langle
\nabla_YX, Z\rangle=-\langle X,
\nabla_YZ\rangle = -\langle X, \nabla_ZY\rangle \\
& = & \langle \nabla_ZX,Y\rangle= \langle
\nabla_XZ,Y\rangle=-\langle Z,\nabla_XY\rangle.\end{eqnarray*} Then
$\langle \nabla_XY, Z\rangle=0$ for any
$Z\in\Gamma(T^{(1,0)}S^{2n})$, this shows $
\nabla_XY\in\Gamma(T^{(1,0)}S^{2n})$.
 \ \ $\Box$

The operator $J$ gives a section $f\colon\; S^{2n}\to
\widetilde{\cal J}(S^{2n})$,
$$f=e_{0}\cdot e_{-1}^t- e_{-1}\cdot e_{0}^t
+\sum\limits_{i=1}^n \ (e_{2i}\cdot e_{2i-1}^t- e_{2i-1}\cdot
e_{2i}^t),$$ and $J=J_f=\sum\limits_{i=1}^n \ (e_{2i}\cdot
e_{2i-1}^t- e_{2i-1}\cdot e_{2i}^t)\colon\; \ TS^{2n}\to TS^{2n}.$
Hence we have
$$f^*\alpha_{ij}^*=\frac 14f^*\langle dA,\alpha_{ij}\rangle
=\frac 14\langle df,\alpha_{ij}|_f\rangle=
\omega_{2i}^{2j-1}+\omega_{2i-1}^{2j},$$
$$f^*\beta_{ij}^*=
\omega_{2i}^{2j}-\omega_{2i-1}^{2j-1}, \ \ \
f^*\omega^{2i-1}=\omega^{2i-1},\ \ \ f^*\omega^{2i}=\omega^{2i}.$$
By Lemma 3.2, we have $Jf^*= f^*\tilde J$ if $J$ is integrable, then
$f\colon\;  S^{2n}\to \widetilde{\cal J}(S^{2n})$ is holomorphic and
$S^{2n}$ is also a Kaehler manifold. This can occur only when $n=1$.
We have proved

{\bf Theorem 3.3} \ There is no orthogonal  complex structure on
$S^{2n}$ for $n>1$.

Now we study the equation $\tilde Jf_*=f_*J$.

Let $T^V\widetilde{\cal J}(S^{2n}), T^H\widetilde{\cal J}(S^{2n})$
be the subspaces of $T\widetilde{\cal J}(S^{2n})$ generated by
$\alpha_{ij},\beta_{ij}$ and $\widetilde X_{2i-1},\widetilde X_{2i}$
respectively. It is easy to see that $T^V\widetilde{\cal J}(S^{2n})$
are tangent to the fibres of fibre bundle $\pi\colon\;
\widetilde{\cal J}(S^{2n})\to S^{2n}$ and
$$\pi_*(\widetilde X_{2i-1})=e_{2i-1}, \ \pi_*(\widetilde X_{2i})=e_{2i}.$$ Note that
$$ \sum\limits_{l=1}^{2n} \ \omega^l\widetilde X_{l} = de_0e_{-1}^t - e_{-1}de_{0}^t
 - Ade_0e_{0}^t + e_{0}de_{0}^tA.$$
The subspaces $T^V\widetilde{\cal J}(S^{2n})$ and
$T^H\widetilde{\cal J}(S^{2n})$  are orthogonal with respect to the
Riemannian metric on $\widetilde{\cal J}(S^{2n})$ and are called the
vertical and horizontal subspaces of $T\widetilde{\cal J}(S^{2n})$
respectively.

By \S2,   any $A\in\widetilde{\cal J}(S^{2n})$ defines a complex
structure $J_A$ on $T_{e_0}S^{2n}$:
$$ X\in T_{e_0}S^{2n}\to  J_A(X)=AX, \ \ e_0=\pi(A).$$

{\bf Lemma 3.4} \ (1) The following diagram of maps is commutative,
$$ \begin{array}{cccccc} T_{A}\widetilde{\cal J}(S^{2n}) & \stackrel{\tilde J}
\longrightarrow  & T_{A}\widetilde{\cal J}(S^{2n}) \\
\pi_{*} \downarrow \quad & {} & \quad \downarrow \pi_{*} \\
T_{e_0}S^{2n} & \stackrel{J_A} \longrightarrow & T_{e_0}S^{2n};
\end{array}$$
 (2)  The complex structure  $\tilde J$   preserves the decomposition
$$T\widetilde{\cal J}(S^{2n})=T^V\widetilde{\cal J}(S^{2n})\oplus T^H\widetilde {\cal
J}(S^{2n}).$$

Let $f\colon\; S^{2n}\to \widetilde{\cal J}(S^{2n})$ be a section
(or a local section) and $J_f$ be the related almost complex
structure on $S^{2n}$. From Lemma 3.4 and $\pi\circ f=id$, we have

{\bf Lemma 3.5} \ $J_f=\pi_*\tilde Jf_*.$

The section $f\colon\; S^{2n}\to \widetilde{\cal J}(S^{2n})$ can be
viewed as a map $f\colon\; \ S^{2n} \to  gl(2n+2,R)$ and
$T_A\widetilde{\cal J}(S^{2n})$ is a subspace of $  gl(2n+2,R)$,
then
$$f_*X=(df)X= Xf, \ \ \ X\in \Gamma (TS^{2n}\otimes C).$$
On the other hand,  $J_f\in \Gamma(End(TS^{2n}))$, we can compute
$\nabla_X J_f$, where $\nabla$ is the Riemannian connection on
$S^{2n}$.

{\bf Lemma 3.6} \ $f_*X=\nabla_X J_f + \widetilde X$, where
$\widetilde X$ is the horizontal lift of $X$.

{\bf Proof}  \ Let $e_1,\cdots,e_{2n}$ be local orthonormal vector
fields on $S^{2n}$. The section $f$ can be represented by
$$f= e_{0}\cdot e_{-1}^t- e_{-1}\cdot e_{0}^t +
(e_1,\cdots,e_{2n})B(e_1,\cdots,e_{2n})^t,$$
$$J_f
=(e_1,\cdots,e_{2n})B(e_1,\cdots,e_{2n})^t, $$
 where $e_0\in S^{2n}$
and $B$  a matrix function on $S^{2n}, \ B^2=-I$ and $BB^t=I$. Let
$$de_0=(e_1,\cdots,e_{2n})(\omega^1 ,\cdots,\omega^{2n})^t,$$
$$d(e_1,\cdots,e_{2n})=(e_1,\cdots,e_{2n})\omega - e_0(\omega^1 ,\cdots,\omega^{2n}).$$
Hence $$\nabla (e_1,\cdots,e_{2n})=(e_1,\cdots,e_{2n})\omega,$$
$$df=(e_1,\cdots,e_{2n})(dB+\omega
B-B\omega)(e_1,\cdots,e_{2n})^t + \sum\limits_{l=1}^{2n} \
\omega^l\widetilde X_{l},$$
$$\nabla J_f = (e_1,\cdots,e_{2n})(dB+\omega B-B\omega)(e_1,\cdots,e_{2n})^t.$$
 These shows
$$f_*X =Xf =\nabla_XJ_f +\widetilde X, $$
where $\widetilde X$ is the horizontal lift of $X$ and $\nabla_XJ_f
\in \Gamma(T^V\widetilde{\cal J}(S^{2n}))$. \ \ \ $\Box$

Now assuming that $J_f$ is integrable. Let
$X,Y\in\Gamma(T^{(1,0)}S^{2n})$,  $\nabla_XY, \ \nabla_YX\in
\Gamma(T^{(1,0)}S^{2n})$  by Lemma 3.2. Then we have
$$(\nabla_XJ_f)Y=\nabla_X(J_fY)-J_f\nabla_XY=0.$$
For any $Y_1\in \Gamma(S^{2n}), \ Y=(1-\sqrt{-1}J_f)Y_1\in
\Gamma(T^{(1,0)}S^{2n})$, we have
$$(\nabla_XJ_f)Y=(\nabla_XJ_f-\sqrt{-1}\nabla_XJ_f\cdot J_f)Y_1=0.$$
Then
$$\nabla_XJ_f-\sqrt{-1}\nabla_XJ_f\cdot J_f=0,$$
$$\nabla_XJ_f=\sqrt{-1}\nabla_XJ_f\cdot J_f=-\sqrt{-1}J_f\nabla_XJ_f,$$
$$\tilde J\nabla_XJ_f=J_f\nabla_XJ_f=\sqrt{-1}\nabla_XJ_f.$$

These shows $f_*X=\nabla_X J_f + \widetilde
X\in\Gamma(T^{(1,0)}\widetilde{\cal J}(S^{2n}))$ for any $X\in
\Gamma(T^{(1,0)}(S^{2n}))$ and the  map $f\colon\; S^{2n}\to
\widetilde{\cal J}(S^{2n})$ is holomorphic. These also proves the
Theorem 3.3.

\vskip 0.7cm \noindent{\bf  4. The complex structure
 on the sphere $S^{2n}$} \vskip
0.2cm

 In the following we use   Riemannian connection  to study the complex structure.
The metric on $S^{2n}$ can be represented by
$$ds^2=\frac {1}{(1+\frac 14 |y|^2)^2} \sum\limits_{i=1}^{2n} \
(dy^i)^2,$$ where $(y^1,y^2,\cdots,y^{2n})$ are the coordinates on
$S^{2n}$ defined by the stereographic projection. Let
$e_i=(1+\frac 14 |y|^2)\frac {\partial}{\partial y^i}$ be an
orthonormal frame fields on $S^{2n}$,\ $\omega^i=\frac 1{1+\frac
14 |y|^2}dy^i$ be their dual 1-forms.

{\bf Lemma 4.1}  \ The Riemannian connection $\nabla$ on $S^{2n}$
is defined by
$$\nabla_{e_i}e_j= -\frac 12 y^je_i +\frac 12 \sum \
y^ke_k\delta_{ij}.$$

{\bf Proof}  \ By the structure equations of Riemannian connection
$d\omega^k=\sum \ \omega^j\wedge\omega_j^k, \
\omega_j^k+\omega_k^j=0,$ we have
$$\omega_j^k=-\frac 12y^j\omega^k + \frac 12y^k\omega^j, \ \
k,j=1,2,\cdots,2n.$$ The Riemannian connection $\nabla$ on
$S^{2n}$ is defined by $\nabla e_j=\sum \ \omega_j^ke_k.$ \ \ \
$\Box$

In the following we sometimes omit the notation $\sum$.

Let $J$ be an almost complex structure on $S^{2n}, \ Je_i=\sum \
e_jB_{ji}$.

{\bf Lemma 4.2}  \   We have
\begin{eqnarray*} &&\nabla_{e_i-\sqrt
{-1}Je_i}(e_j-\sqrt {-1}Je_j) \\
&& =-\sqrt {-1}[(e_i-\sqrt {-1}Je_i)B_{kj}]e_k-
 \frac 12(y^j-\sqrt{-1}y^kB_{kj})(e_i-\sqrt {-1}Je_i) \\
&& \quad +\frac 12[\delta_{ij} -B_{ki}B_{kj} -
\sqrt{-1}(B_{ij}+B_{ji})] y^le_l. \end{eqnarray*}

Then  \begin{eqnarray*} && [e_i-\sqrt {-1}Je_i,e_j-\sqrt {-1}Je_j] \\
&& =-\sqrt {-1} [(e_i-\sqrt {-1}Je_i)B_{kj}]e_k+\sqrt {-1}
[(e_j-\sqrt {-1}Je_j)B_{ki}]e_k \\
&& \quad -
 \frac 12(y^j-\sqrt{-1}y^kB_{kj})(e_i-\sqrt {-1}Je_i)+
 \frac 12(y^i-\sqrt{-1}y^kB_{ki})(e_j-\sqrt {-1}Je_j).\end{eqnarray*}

{\bf Corollary 4.3}  \ (1)  $J$ is integrable if and only if
$$\sum \ [e_iB_{kj}- e_jB_{ki}]e_k+
\sum \ [(Je_i)B_{kj}- (Je_j)B_{ki}]Je_k=0;$$

(2) If $J$ is orthogonal, $J$ is  integrable if and only if each
$\sum \ [(e_i-\sqrt {-1}Je_i)B_{kj}]e_k$ is $(1,0)$, then
$$\sum \ (e_iB_{kj})e_k+\sum \ (Je_i)B_{kj}Je_k  =0.$$

{\bf Proof}  \ (1) By Lemma 4.2, $[e_i-\sqrt {-1}Je_i,e_j-\sqrt
{-1}Je_j]$ is $(1,0)$ if and only if $$\sum \ [(e_i-\sqrt
{-1}Je_i)B_{kj}]e_k-\sum \ [(e_j-\sqrt {-1}Je_j)B_{ki}]e_k $$ is
$(1,0)$. Act $J$ on the   equation we have
$$\sum \ [(e_iB_{kj})e_k- (e_jB_{ki})e_k]=
-\sum \ [((Je_i)B_{kj})Je_k- ((Je_j)B_{ki})Je_k].$$

(2) follows from Lemma 3.2 and 4.2. \ \ \ $\Box$

{\bf Remark}  \ If $J$ is orthogonal and integrable, we have
$BB^t=I, \ B=-B^t$, then
\begin{eqnarray*} &&\nabla_{e_i-\sqrt
{-1}Je_i}(e_j-\sqrt {-1}Je_j) \\
&& =-\sqrt {-1}[(e_i-\sqrt {-1}Je_i)B_{kj}]e_k-
 \frac 12(y^j-\sqrt{-1}y^kB_{kj})(e_i-\sqrt {-1}Je_i) \\
&& =-(Je_i)B_{kj}(e_k-\sqrt {-1}Je_k)-
 \frac 12(y^j-\sqrt{-1}y^kB_{kj})(e_i-\sqrt {-1}Je_i) \\
&& \not= -\nabla_{e_j-\sqrt {-1}Je_j}(e_i-\sqrt {-1}Je_i).
\end{eqnarray*}
  But   Lebrun  [4]  proved
$\nabla_{e_i+\sqrt {-1}Je_i}(e_j+\sqrt {-1}Je_j) = -
\nabla_{e_j+\sqrt {-1}Je_j}(e_i+\sqrt {-1}Je_i)$.

Let $f=\tilde J= e_{0}\cdot e_{-1}^t- e_{-1}\cdot e_{0}^t +
(e_1,\cdots,e_{2n})B(e_1,\cdots,e_{2n})^t$ be a section of twistor
space ${\cal J}(S^{2n})$, $B^2=-I$. We also view $f$ as  a map from
$S^{2n}$ to the Kaehler manifold ${\cal J}(R^{2n+2})$. Then
$J=J_f=\sum \ e_iB_{ij}e_j^t$ is an almost complex structure on
$S^{2n}$. By \S 3, we have
$$f_*e_l=e_lf =\nabla_{e_l}J_f +\widehat e_l,$$
where
\begin{eqnarray*}  \nabla_{e_l}J_f & = &
 e_i(e_lB_{ij})e_j^t + (\nabla _{e_l}e_i)B_{ij}e_j^t+
e_iB_{ij}(\nabla _{e_l}e_j)^t \\
&=&  e_i(e_lB_{ij})e_j^t + \frac 12(-e_ly^ie_i^t +y^ke_ke_l^t)
J+\frac 12 J(e_ly^ie_i^t -y^ie_ie_l^t),\end{eqnarray*}
$$ \widehat e_l=e_le_{-1}^t-e_{-1}e_l^t+ \tilde J
(e_le_{-1}^t-e_{-1}e_l^t) \tilde J.$$

{\bf Theorem 4.4}  \ Let $f\colon\; S^{2n}\to {\cal
J}(S^{2n})\subset {\cal J}(R^{2n+2})$ be a local section. The map
$f$ is holomorphic if and only if $f$ is a local section of
$\widetilde{\cal J}(S^{2n})$  and integrable.

{\bf Proof}  \ For any tangent vector $e_l$ defined above, we have
\begin{eqnarray*}  (\tilde J f_*-f_*J_f)e_l & = &
Je_i(e_lB_{ij})e_j^t - e_i((Je_l)B_{ij})e_j^t \\
&& -\tilde J (e_0-\frac 12y^ie_i)(B_{lj}+B_{jl})e_j^t
 +(e_0- \frac 12y^ie_i)(B_{jl}+B_{lj})e_j^tJ. \end{eqnarray*}
This shows the local section $f$ is holomorphic if and only if
$B_{lj}+B_{jl}=0$ and $Je_i(e_lB_{ij}) - e_i((Je_l)B_{ij})=0$. That
is, the complex structure $J_f$ is orthogonal and integrable. \ \ \
$\Box$

If the complex structure $J_f$ is orthogonal and integrable,
$B_{ij}+B_{ji}=0$, we have
\begin{eqnarray*} && \nabla_{e_l-\sqrt {-1}Je_l}J_f \\
&&= e_i[(e_l-\sqrt {-1}Je_l)B_{ij}]e_j^t +\frac {\sqrt
{-1}}2(1-\sqrt {-1} J)e_l y^je_j^t(1+\sqrt {-1}
J) \\
&&  \quad - \frac {\sqrt{-1}}2 y^i(1-\sqrt {-1} J)e_ie_l^t(1
+\sqrt{-1}J),\end{eqnarray*} and \begin{eqnarray*} \widehat
{e_l-\sqrt {-1}Je_l} &=& (1-\sqrt {-1}J)e_le_{-1}^t(1+\sqrt
{-1}\tilde J) \\ && - (1-\sqrt {-1}\tilde J)e_{-1}(e_{l}^t-\sqrt
{-1}B_{kl}e_{k}^t).
\end{eqnarray*}
 By Corollary 4.3(2),  $e_i[(e_l-\sqrt
{-1}Je_l)B_{ij}]e_j^t$ is $(1,0)$ vector field, then $f_*(e_l-\sqrt
{-1}Je_l)$ is $(1,0)$  vector field on $\widetilde{\cal J}(S^{2n})$.
These gives another proof of Theorem 3.3.

\vskip 0.7cm \noindent{\bf  5. The map $f: S^{2n}\to G_C(n, 2n+1)$}
\vskip 0.2cm

The complex Grassmann manifold $G_C(k, n)$ is a Kaehler manifold
with   complex structure $\tilde J$. Let
$s_1,\cdots,s_k,s_{k+1},\cdots, s_n$ be complex functions on
 an open subset $U$ of $G_C(k, n)$ with values in
$C^n$ such that the elements of $U$ can be generated  by
$s_1,\cdots,s_k$, and $s_1, \cdots, s_n$ are   linearly independent
everywhere. Let $$ds_i=\sum\limits_{B=1}^n \ \widetilde
\varphi_i^Bs_B, \ i=1,\cdots,k.$$
 We have
$$d(s_1\wedge \cdots\wedge s_k)=\sum\limits_{i=1}^k \ \widetilde
\varphi_i^is_1\wedge \cdots\wedge s_k
+\sum\limits_{i=1}^k\sum\limits_{\alpha=k+1}^n \ \widetilde
\varphi_i^\alpha  s_1\wedge\cdots s_{i-1}\wedge s_\alpha \wedge
s_{i+1}\cdots\wedge s_k.$$

{\bf Lemma 5.1}  \ $\{\widetilde \varphi_i^\alpha\}$ are basis of
$T^{(1,0)}G_C(k, n)$ on $U$.

{\bf Proof} \ Let $$s_i=(0,\cdots,\underset{i}1,\cdots,0,z_{i
k+1},\cdots, z_{i n}), \ \ i=1,\cdots,k, \ z_{i\alpha}\in C,$$
$$s_\alpha=(0,\cdots,0,\underset {\alpha}1,0,\cdots,0), \ \alpha=k+1,\cdots,n.$$
 $s_1,\cdots,s_k$ generate an open subset $U$ of
$G_C(k, n)$   with the complex coordinates $\{z_{i\alpha}\}$. In
these case, $\widetilde \varphi_i^\alpha=dz_{i\alpha}$ form a basis
of $T^{(1,0)}G_C(k, n)$ at every point of $U$.  The general cases
can be proved easily. See also [11].\ \ \ $\Box$

 Let $J$ be an orthogonal almost complex structure on
$S^{2n}$, $e_0, e_1,\cdots, e_{2n}$ and $Z_i=e_{2i-1}-\sqrt
{-1}e_{2i}, \ Z_{\bar i}=e_{2i-1}+\sqrt {-1}e_{2i}$ are defined as
in \S 3.

As [4], define  map $f: S^{2n}\to G_C(n, 2n+1)$,
 $f(p)=T_p^{(0,1)}S^{2n},  \ p\in S^{2n}$.
Let $s_1, \cdots, s_{2n+1}$ be local  frame fields on  $G_{\bf C}(n,
2n+1)$ such that restrict on $f(S^{2n})$ we have
 $s_i=Z_{\bar i}, \ s_{n+i} = Z_i, \  i=1,\cdots,n, \ s_{2n+1}=e_0$. Hence
$$f^*(\widetilde \varphi_i^{n+j})=\varphi_{\bar i}^j, \
f^*(\widetilde \varphi_i^{2n+1})=-\varphi^i, \ i,j=1,\cdots,n.$$
Then by Lemma 3.2 and  Lemma 5.1, we have $f^*\widetilde J=Jf^*$ on
$T^{(1,0)}G_{\bf C}(n, 2n+1)$ if $J$ is integrable, hence the map
$f: S^{2n}\to G_{\bf C}(n, 2n+1)$ is holomorphic. This also proves
Theorem 3.3.

Note that $[s_1\wedge \cdots\wedge s_k]$ give an imbedding of
$G_C(k, n)$ in complex projective space ${\bf P}(\bigwedge^k(C^n))$.
$G_C(k, n)$ is a complex submanifold of  ${\bf
P}(\bigwedge^k(C^n))$, we can also use the map $f: S^{2n}\to {\bf
P}(\bigwedge^n(C^{2n+1}))$ to study the problem.

\newpage

\vskip 0.7cm \noindent{\bf  6. The first Chern class of
$T^{(1,0)}S^{2n}$} \vskip 0.2cm

In \S 2, we   show that  the orthogonal twistor space
$\widetilde{\cal J}(S^{2n})$ of the sphere $S^{2n}$ is a Kaehler
manifold. We can show that the first  Chern class of the   vector
bundle $T^{H(1,0)} \widetilde{\cal J}(S^{2n})$  can be represented
by the Kaehler form of $ \widetilde{\cal J}(S^{2n})$ with a constant
as  coefficient. Thus $c_1(T^{(1,0)}S^{2n})=f^*c_1(T^{H(1,0)}
\widetilde{\cal J}(S^{2n}))\in H^2(S^{2n})$ is non-zero if the map
$f\colon\; S^{2n}\to \widetilde{\cal J}(S^{2n}) $ is holomorphic. In
 following we use $c_1(T^{(1,0)}S^{2n})$ to  study the complex
structure on sphere $S^{2n}$.

As  \S 3, let  $J$ be an orthogonal almost complex structure on
$S^{2n}$, $e_1,\cdots, e_{2n}$ be local orthonormal frame fields
with $Je_{2i-1}=e_{2i},\ i=1,\cdots,n$. The  Riemannian connection
on $S^{2n}$ is $\nabla e_B=\sum\limits_{C=1}^{2n} \ \omega_B^Ce_C$,
then $$\nabla^2 e_B=\sum\limits_{C=1}^{2n} \
\Omega_B^Ce_C=\sum\limits_{C=1}^{2n} \
(d\omega_B^C-\sum\limits_{D=1}^{2n} \
\omega_B^D\wedge\omega_D^C)e_C,$$ where $\omega_B^C$ is defined in
\S 3. The connection on $T^{(1,0)}S^{2n}$ is
$$\nabla
Z_{i}=\sum\limits_{j=1}^n \ \varphi_{i}^{j} Z_{j}, \ \ \nabla^2
Z_{i}=\sum\limits_{k=1}^n \ \Phi_i^kZ_k,$$ where $
\varphi_{i}^{j}=\frac 12[\omega_{2j-1}^{2i-1}+
\omega_{2j}^{2i}+\sqrt{-1}(\omega_{2j-1}^{2i}-\omega_{2j}^{2i-1})],
\ \Phi_i^k=\mbox{d} \varphi_i^k-\sum \
\varphi_i^j\wedge\varphi_j^k$.

As $\sum\limits_{i,j} \
\varphi_i^j\wedge\varphi_j^i=-\sum\limits_{i,j} \
\varphi^i_j\wedge\varphi^j_i=0$, we have
$$c_1(T^{(1,0)}S^{2n})=\frac {\sqrt{-1}}{2\pi}\sum\limits_{i=1}^n \ \Phi_i^i=
\frac {\sqrt{-1}}{2\pi}\sum\limits_{i=1}^n \ d\varphi_i^i=-\frac
{1}{2\pi}\sum\limits_{i=1}^n \ d\omega_{2i-1}^{2i}.$$

Let   $A=(\omega_{2i-1}^{2j-1}),\ B=(\omega_{2i-1}^{2j}),\
C=(\omega_{2i}^{2j-1}),\ D=(\omega_{2i}^{2j}), \ \omega =\left(
\begin{array}{cccccccc}
A  & B  \\
C & D
\end{array}\right), \ \Omega =d\omega-\omega\wedge \omega$ and
$J_0=\left(
\begin{array}{cccccccc}
{}  & -I  \\
I & {}
\end{array}\right)$. By $$ \sum\limits_{i=1}^n \ \mbox{d} \omega_{2i-1}^{2i}= \frac 12\mbox{tr}(\mbox{d}
\omega J_0)=\frac {1}{2}\mbox{tr} ( \Omega J_0+  \omega \wedge\omega
J_0),$$ we have

{\bf  Lemma 6.1} \  The first Chern class of $T^{(1,0)}S^{2n}$ can
be represented by
$$c_1(T^{(1,0)}S^{2n})
=-\frac {1}{4\pi}\mbox{tr} ( \Omega J_0+  \omega \wedge\omega
J_0).$$

{\bf  Theorem 6.2} \ When $n>1$, there is no complex structure in a
neighborhood of the space $ \widetilde{\cal J}(S^{2n})$.

{\bf Proof} \  Assuming  $J$ is integrable,   $J(B+C)=A-D$ by Lemma
3.2. By
 $$\omega +J_0 \omega J_0=\left(
\begin{array}{cccccccc}
A-D  & B+C  \\
B+C & -A+D
\end{array}\right),\ J_0 \omega - \omega J_0=\left( \begin{array}{cccccccc}
-B-C  & A-D  \\
A-D & B+C
\end{array}\right),$$ we have $J( \omega +J_0  \omega
J_0)=J_0(\omega +J_0 \omega J_0).$ The sectional curvature of the
metric  on $S^{2n}$ is constant,  $ \Omega_k^l=- \omega^k\wedge
\omega^l.$ Then  we have
\begin{eqnarray*} && \mbox{tr}( \Omega J_0)+\mbox{tr}( \omega\wedge \omega J_0) \\
&&=\mbox{tr}( \Omega J_0)+\frac 14\mbox{tr}[( \omega +J_0 \omega J_0)\wedge( \omega J_0-J_0 \omega)] \\
&&= \mbox{tr}( \Omega J_0)-\mbox{tr}[(B+C)\wedge J(B+C)] \\
&&= \mbox{tr}( \Omega J_0)+\mbox{tr}[(B+C)\wedge J(B+C)^t] \\
 &&=2\sum\limits_{i=1}^n \  \omega^{2i-1}\wedge J\omega^{2i-1}
 + 2\sum\limits_{i<j} \ ( \omega_{2i-1}^{2j}+ \omega_{2i}^{2j-1})\wedge J
( \omega_{2i-1}^{2j}+ \omega_{2i}^{2j-1}).
\end{eqnarray*}

For any $X \in TS^{2n}$, we have
\begin{eqnarray*}
&& [\mbox{tr}( \Omega J_0)+\mbox{tr}( \omega\wedge \omega J_0)](X, JX)   \\
 &  &=-2\sum\limits_{i=1}^n \  ([\omega^{2i-1}(X)]^2+[\omega^{2i}(X)]^2) \\
 && \quad
- 2\sum\limits_{i<j} \ [( \omega_{2i-1}^{2j}+
\omega_{2i}^{2j-1})(X)]^2 - 2\sum\limits_{i<j} \ [(
\omega_{2i-1}^{2j-1}- \omega_{2i}^{2j})(X)]^2.
\end{eqnarray*}
Then  2-form $\mbox{tr}(\Omega J_0+ \omega\wedge \omega J_0)$ are
non-degenerate everywhere and $S^{2n}$ becomes a symplectic
manifold, this contradict to the fact of $H^2(S^{2n})=0$ for $n>1$.
As the Riemannian curvature is continuous with  the Riemannian
metric, these shows there is a neighborhood of $\widetilde{\cal
J}(S^{2n})$ in ${\cal J}(S^{2n})$ such that there is no complex
structure in this neighborhood. \ \ \ \ \ $\Box$

\vskip 1cm \noindent{\large \bf References} \vskip 0.2cm

{\small

\noindent [1] \  Albuquerque, R.,  Salavessa, I. M. C.: On
the twistor space of pseudo-spheres.
Differential Geomtry and its applications, {\bf 25}, 207-221(2007)

\noindent [2] \  Borel,   A.,   Serre, J. P.: Groupes de Lie et
puissances r\'{e}duites de Steenrod. Amer. J. Math. {\bf 75},
409-448(1953)

\noindent [3] \ Lawson, Jr  H. B., Michelsohn,  M.:  Spin geometry.
Princeton University Press. Princeton New Jersey, 1989

\noindent [4] \  Lebrun, C.: Orthogonal complex structure on $S^6$.
Proc. of the Amer. Math. Soc., {\bf 101}, 136-138(1987)

\noindent [5] \  O'Brian, N. R.,  Rawnsley, J.: Twistor spaces. Ann.
of Global Analysis and Geometry. {\bf 3}, 29-58(1985)

\noindent [6] \ Steenrod, N.: The topology of fibre bundles,
Princeton Univ. Press, 1974.

 \noindent [7] \  Zhou, J. W.: Grassmann manifold $G(2,8)$
and Complex Structure on $S^{6}$.  arXiv: math. DG/0608052

\noindent [8] \  Zhou, J. W.: Morse functions on Grassmann
manifolds. Proc. of the Royal Soc. of Edinburgh, {\bf 135A}(2005),
209-221

\noindent [9] \  Zhou, J. W.: The geometry and topology on Grassmann
manifolds.  Math. J. Okayama Univ., {\bf 48}(2006), 181-195

\noindent [10] \  Zhou, J. W.: A note on characteristic classes.
Czechoslovak Math. J., {\bf 58}(2006), 721-732  }

 \noindent [11] \  Zhou, J. W.: Lectures on differential geometry (in Chinese).
 Science Press, Beijing, 2010.

\end{document}